\documentclass[12pt]{amsart}

\usepackage[left=2cm, right=2cm, top=3cm, bottom=3cm]{geometry}

\usepackage{amsmath, amsthm, amsfonts, amssymb}
\usepackage{mathrsfs}
\usepackage{mathtools}
\usepackage{upgreek}
\usepackage{booktabs}
\usepackage{url}

\usepackage{listings}
\makeatletter
\@namedef{subjclassname@2020}{%
  \textup{2020} Mathematics Subject Classification}
\makeatother

\theoremstyle{plain}
\newtheorem{thm}{Theorem}[section]
\newtheorem*{thm*}{Main Theorem}{}
\newtheorem{coro}[thm]{Corollary}
\newtheorem{prop}[thm]{Proposition}

\newtheorem{lemm}[thm]{Lemma}

\theoremstyle{definition}

\theoremstyle{remark}
\newtheorem{rema}[thm]{Remark}

\newtheorem{conj}[thm]{Conjecture}

\newcommand\twoscript[2]{\substack{{#1} \\ {#2}}}

\newcommand\rmi{\mathrm{i}}

\newcommand\tbtmat[4]{\left(\begin{smallmatrix}{#1} & {#2} \\ {#3} & {#4}\end{smallmatrix}\right)}
\newcommand\tbtMat[4]{\begin{pmatrix}{#1} & {#2} \\ {#3} & {#4}\end{pmatrix}}
\newcommand*\abs[1]{\lvert#1\rvert}
\newcommand\etp[1]{\mathfrak{e}\left(#1\right)}

\newcommand\numZ{\mathbb{Z}}

\newcommand\numR{\mathbb{R}}
\newcommand\numC{\mathbb{C}}

\newcommand\numgeq[2]{\mathbb{#1}_{\geq #2}}

\newcommand\slZ{\mathrm{SL}_2(\mathbb{Z})}

\newcommand\pslZ{\mathrm{PSL}_2(\mathbb{Z})}

\newcommand\sltZ{\widetilde{\mathrm{SL}_2(\mathbb{Z})}}

\newcommand\uhp{\mathfrak{H}}

\newcommand\elesltRaDs[6]{\left(\left(\begin{smallmatrix}{#1} & {#2} \\ {#3} & {#4}\end{smallmatrix}\right),{#5}\right)}
\newcommand\Dcover[2]{\widetilde{#1^{#2}}}

\begin{document}


\baselineskip=17pt


\title[Linear characters of $\Gamma_0(N)$]{Explicit formulae for linear characters of $\Gamma_0(N)$}

\author{Xiao-Jie Zhu}
\address{School of Mathematical Sciences\\
East China Normal University\\
500 Dongchuan Road, 200241\\
Shanghai, P. R. China}
\email{zhuxiaojiemath@outlook.com}
\urladdr{https://orcid.org/0000-0002-6733-0755}

\date{}

\begin{abstract}
We give explicit formulae for a class of complex linear unitary characters of the congruence subgroups $\Gamma_0(N)$ which involve a variant of Rademacher's $\Psi$ function. We then prove that these characters cover all characters of $\Gamma_0(N)$ precisely when $N=1,2,3,4,5,6,7,8,10,12,13$.
\end{abstract}



\thanks{This work is supported by Key Laboratory of Mathematics and Engineering Applications, Ministry of Education, P. R. China and by Shanghai Key Laboratory of Pure Mathematics and Mathematical Practice.}

\maketitle


\section{Introduction}
\label{sec:Introduction}
The congruence subgroup
\begin{equation*}
\Gamma_0(N)=\left\{\tbtmat{a}{b}{c}{d}\in\slZ\colon c \equiv 0 \bmod N\right\}
\end{equation*}
plays an important role in number theory, especially in the theory of modular forms, elliptic curves, and quadratic forms. The structures, more precisely, presentations, representations and free-product decompositions of these groups or their projections onto transformation groups have been studied for a long while. Rademacher \cite{Rad29} gave a presentation of $\Gamma_0(N)$ when $N$ is a prime and Chuman \cite{Chu73} generalized Rademacher's presentation to any integer $N$. Unfortunately Chuman omitted one relator in his presentation. Lascurain \cite{Las02} fixed this mistake and further simplified Chuman's presentation. Their results are all based on Reidemeister-Schreier method (c.f. \cite[Proposition 4.1, Chap II]{LS77}).

Rademacher also noted that $\pslZ=\slZ/\{\pm I\}$ is (isomorphic to) a free product of $\numZ/2\numZ$ and $\numZ/3\numZ$; hence by Kurosh subgroup theorem (c.f. \cite[eq. (3)]{Mac58}) any of its finite-index subgroup is a free product of finitely many copies of $\numZ/2\numZ$, $\numZ/3\numZ$ and $\numZ$. Kulkarni \cite{Kul91}, in 1991, introduced the concept of generalized Farey sequences and Farey symbols and obtained an algorithm for factoring any $\Gamma_0(N)/\{\pm I\}$ into a free product. His method is based on fundamental domains, side-pairing transformations and hyperbolic geometry, which is completely different from Lascurain's.

Representations also contain information on the structure of $\Gamma_0(N)$. In this paper, we shall consider one-dimensional complex unitary representations of $\Gamma_0(N)$. Throughout the whole paper, by the term \emph{character} we mean a one-dimensional complex unitary representation. Well-known examples are characters induced by Dirichlet characters modulo $N$ and characters of some eta-quotients of level $N$ (c.f. \cite[eq. (15) and (16)]{ZZ23}), which are widely used in the theory of modular forms. However, these are not all characters. The aim of this paper is to give some characters of $\Gamma_0(N)$ using explicit and elementary expressions and to show that these characters form the whole group of characters precisely when $N=1,2,3,4,5,6,7,8,10,12,13$.

Our formulae involve Dedekind sums, which are defined by
\begin{equation*}
s(h,k)=\sum_{r=1}^{k-1}\frac{r}{k}\left(\frac{hr}{k}-\left[\frac{hr}{k}\right]-\frac{1}{2}\right),
\end{equation*}
where $h$ and $k$ are coprime integers and $k > 0$, and $[x]$ means the greatest integer smaller than or equal to $x$. Dedekind sums first occurred in the study of transformation equations of Dedekind $\eta$ function and possess many good properties and identities (c.f. \cite{RM72}). In 1956, Rademacher \cite{Rad56} introduced a function which he called $\Psi$ and proved interesting formulae for Dedekind sums using this function. For our purpose, we need a slightly different version\footnote{Let Rademacher's function \cite[eq. (11)]{Rad56} be denoted by $\Psi_R$. Then for $\tbtmat{a}{b}{c}{d}\in\slZ$, $\Psi\tbtmat{a}{b}{c}{d}=\Psi_R\tbtmat{a}{b}{c}{d}$ if and only if $a+d>0$. For $a+d<0$ ($=0$ resp.) we have $\abs{\Psi\tbtmat{a}{b}{c}{d}-\Psi_R\tbtmat{a}{b}{c}{d}}=6$ ($=3$ resp.).}, which we also call $\Psi$ as follows:
\begin{align}
\label{eq:Psi}
\Psi\colon \slZ &\rightarrow \numZ \\
\tbtmat{a}{b}{c}{d} &\mapsto \begin{dcases}
\frac{a+d}{c}+12s(-d,c)-3,   & c>0; \\
\frac{a+d}{c}+12s(d,-c)+3,   & c<0; \\
b,                           & c=0, a>0; \\
-b-6,                        & c=0, a<0. 
\end{dcases}\notag
\end{align}
The fact that the image of $\Psi$ is contained in $\numZ$ can be proved by taking $24$th powers of both sides of \eqref{eq:etaTrans}. Our main result is the following theorem.
\begin{thm*}
Let $N$ be a positive integer. Then the map
\begin{align}
\label{eq:Gamma0Nchar}
\Gamma_0(N) &\rightarrow S^{1} \\
\tbtmat{a}{b}{c}{d} &\mapsto \chi(d)\etp{r_1\Psi\tbtMat{a}{b}{c}{d}+\sum_{1<l\mid N}r_l\left(\Psi\tbtMat{a}{b}{c}{d}-\Psi\tbtMat{a}{bl}{c/l}{d}\right)}\notag
\end{align}
is a linear unitary character on $\Gamma_0(N)$, where $\chi$ is a Dirichlet character modulo $N$, $r_1\in\frac{1}{12}\numZ$ and $r_l\in\numR$ if $1<l\mid N$. Moreover, all linear unitary characters of $\Gamma_0(N)$ are of this form precisely when $N=1,2,3,4,5,6,7,8,10,12,13$.
\end{thm*}
We explain some notations. The symbol $\etp{x}$ denotes $\exp(2\pi\rmi x)$ and $S^1$ denotes $\{z\in\numC\colon \abs{z}=1\}$. Moreover, we set $T=\tbtmat{1}{1}{0}{1}$, $S=\tbtmat{0}{-1}{1}{0}$, $I=\tbtmat{1}{0}{0}{1}$ and $B_l=\tbtmat{1}{0}{0}{l}$.

It should be noted that, once given a fixed $\Gamma_0(N)$, one may work out all its characters by a presentation, for instance, by the free-product decomposition according to Kulkarni \cite{Kul91}. However, such a character is determined by its values on a set of generators, hence is not explicit. As comparison, one may evaluate the value of \eqref{eq:Gamma0Nchar} at any $\gamma\in\Gamma_0(N)$ immediately, and this does not require a decomposition of $\gamma$ into a product of generators first.

The characters of more general modular groups, for instance, $\mathrm{SL}_2$ over a Dedekind domain, have been determined by Boylan and Skoruppa \cite{BS2013}.

The structure of the paper is as follows. In Section \ref{sec:PsiFunction} we study basic properties of $\Psi$ (c.f. \eqref{eq:Psi}). We then construct homomorphisms from $\Gamma_0(N)$ to $\numZ$ based on $\Psi$ in Section \ref{sec:HomGamma0N}. In Section \ref{sec:groupPresentation}, we recall necessary preliminaries on presentations of Fuchsian groups, especially the theorem of Hoare-Karrass-Solitar \cite[Theorem 3]{HKS71}. These two sections contain key points in stating and proving the Main Theorem, which is done in Section \ref{sec:CharGamma0N}. In the final section, we formulate some related results, open questions and conjectures concerning the main topic.

\section{A variant of Rademacher's $\Psi$ function}
\label{sec:PsiFunction}
The major property of $\Psi$ (c.f. \eqref{eq:Psi}) which we need is the following one.
\begin{prop}
\label{prop:PsiMultiplication}
Let $\gamma_1=\tbtmat{a_1}{b_1}{c_1}{d_1}$ and $\gamma_2=\tbtmat{a_2}{b_2}{c_2}{d_2}$ be matrices in $\slZ$. Set $\gamma_1\gamma_2=\tbtmat{a_3}{b_3}{c_3}{d_3}$. Then
\begin{enumerate}
  \item[(a)] $\Psi(\gamma_1\gamma_2)=\Psi(\gamma_1)+\Psi(\gamma_2)+12$ if one of the following two conditions holds:
  \begin{itemize}
    \item $c_1=c_2=0$ and $d_1<0,\,d_2<0$,
    \item $c_1 \geq 0,\, c_2 \geq 0$ but $c_3 < 0$.
  \end{itemize}
  \item[(b)] $\Psi(\gamma_1\gamma_2)=\Psi(\gamma_1)+\Psi(\gamma_2)-12$ if $c_1 < 0$, $c_2 < 0$ but $c_3 \geq 0$.
  \item[(c)] $\Psi(\gamma_1\gamma_2)=\Psi(\gamma_1)+\Psi(\gamma_2)$ in all other cases.
\end{enumerate}
\end{prop}

Note that Rademacher's original $\Psi$ function also satisfies such properties for special matrices, c.f. \cite[eq. (33)]{Rad56}. The motivation of our modified version is to make such formulae hold for all matrices in $\slZ$. Also note that, one can not modify $\Psi$ to make it a nontrivial homomorphism from $\slZ$ to the additive group $\numC$ according to the well-known presentation of $\slZ$.

To prove Proposition \ref{prop:PsiMultiplication} without tedious verification by cases, we need the concept of a multiple cover of the modular group, slash operators of rational weights, and rational powers of Dedekind $\eta$ function. The reader may refer to \cite[\S 2 and Lemma 4.1]{Zhu23}. We review some elements here.

Let $D$ be a positive integer. The $D$-cover of the full modular group $\Dcover{\slZ}{D}$ is the group consisting of pairs $\elesltRaDs{a}{b}{c}{d}{\varepsilon}{D}$ where $\tbtmat{a}{b}{c}{d} \in \slZ$ and $\varepsilon\in\numC$ satisfies $\varepsilon^D=1$. The law of composition is given by the formula
\begin{equation*}
\left(\tbtmat{a_1}{b_1}{c_1}{d_1}, \varepsilon_1\right)\left(\tbtmat{a_2}{b_2}{c_2}{d_2}, \varepsilon_2\right)=\left(\tbtmat{a_1}{b_1}{c_1}{d_1}\tbtmat{a_2}{b_2}{c_2}{d_2}, \varepsilon_1\varepsilon_2\sigma(\gamma_1,\gamma_2)\right),
\end{equation*}
where $\gamma_i=\tbtmat{a_i}{b_i}{c_i}{d_i}$ and
\begin{equation*}
\sigma(\gamma_1,\gamma_2)=\frac{\sqrt[D]{c_1(a_2\tau+b_2)/(c_2\tau+d_2)+d_1}\sqrt[D]{c_2\tau+d_2}}{\sqrt[D]{c_1(a_2\tau+b_2)+d_1(c_2\tau+d_2)}},\quad \tau \in \numC,\, \Im\tau>0.
\end{equation*}
For the $D$th root of a complex-valued function, we choose the principal branch, so that $\sqrt[D]z=\exp(\frac{1}{D}\log z)$ with $-\uppi < \Im(\log z) \leq \uppi$. One may verify immediately that the cocycle $\sigma(\gamma_1,\gamma_2)$ is independent of $\tau$ and can be evaluated by the following formulae. (Set $\gamma_1\gamma_2=\tbtmat{a_3}{b_3}{c_3}{d_3}$ as in Proposition \ref{prop:PsiMultiplication}.)
\begin{enumerate}
  \item[(a)] $\sigma(\gamma_1,\gamma_2)=\etp{1/D}$ if one of the following two conditions holds:
  \begin{itemize}
    \item $c_1=c_2=0$ and $d_1<0,\,d_2<0$,
    \item $c_1 \geq 0,\, c_2 \geq 0$ but $c_3 < 0$.
  \end{itemize}
  \item[(b)] $\sigma(\gamma_1,\gamma_2)=\etp{-1/D}$ if $c_1 < 0$, $c_2 < 0$ but $c_3 \geq 0$.
  \item[(c)] $\sigma(\gamma_1,\gamma_2)=1$ in other cases.
\end{enumerate}

It seems that the cocycle $\sigma(\gamma_1,\gamma_2)$ can be expressed by Hilbert symbols, c.f. \cite[Theorem 1]{Kub67}. (Kubota's cocycle may be different from ours.) For $D=2$, Str\"{o}mberg gave an expression involving Hilbert symbols in \cite[Theorem 4.1]{Str13}. Note that $\Dcover{\slZ}{2}$ is usually denoted by $\sltZ$ or $\mathrm{Mp}_2(\numZ)$ and is called a metaplectic group.

Let $\uhp$ be the complex upper half-plane, let $k\in \frac{1}{D}\numZ$ and let $f$ be a complex-valued function defined on $\uhp$. The group $\Dcover{\slZ}{D}$ acts on $f$ by $f\vert_k\left(\tbtmat{a}{b}{c}{d}, \varepsilon\right)(\tau)=\varepsilon^{-Dk}(c\tau+d)^{-k}f\left(\frac{a\tau+b}{c\tau+d}\right)$. It follows immediately that $f\vert_k(\gamma_1\gamma_2)=f\vert_k\gamma_1\vert_k\gamma_2$ with $\gamma_1,\gamma_2\in\Dcover{\slZ}{D}$ and $f\vert_k(I,1)=f$. Moreover, if $f$ is not identically zero and there exists $c(\gamma)\in\numC$ such that $f\vert_k\gamma=c(\gamma)f$ for all $\gamma\in\Dcover{\slZ}{D}$, then $\gamma\mapsto c(\gamma)$ must be a character of $\Dcover{\slZ}{D}$.

The well-known transformation equation of the Dedekind $\eta$ function (c.f. \cite[Theorem 3.4]{Apo76}) can be reformulated as (taking $D=2$)
\begin{equation}
\label{eq:etaTrans}
\eta\vert_{1/2}(\gamma,\varepsilon)=\varepsilon^{-1}\cdot\etp{\frac{1}{24}\Psi(\gamma)}\eta,\quad (\gamma,\varepsilon)\in \Dcover{\slZ}{D}.
\end{equation}
For any positive integer $D$, choose a branch of the fractional power $\eta^{2/D}$ as in \cite[Section 4]{Zhu23}. Then
\begin{equation}
\label{eq:etaPowerTrans}
\eta^{2/D}\vert_{1/D}(\gamma,\varepsilon)=\varepsilon^{-1}\cdot\etp{\frac{1}{12D}\Psi(\gamma)}\eta^{2/D},\quad (\gamma,\varepsilon)\in \Dcover{\slZ}{D}.
\end{equation}
The case $2\mid D$ was proved in \cite[Lemma 4.1]{Zhu23}, and the other case can be proved similarly.

\begin{lemm}
Let $D$ be any positive integer, $\gamma_1,\gamma_2\in\slZ$. Then we have
\begin{equation*}
\sigma(\gamma_1,\gamma_2)=\etp{\frac{1}{12D}\left(\Psi(\gamma_1\gamma_2)-\Psi(\gamma_1)-\Psi(\gamma_2)\right)}.
\end{equation*}
\end{lemm}
\begin{proof}
We have $(\eta^{2/D}\vert_{1/D}(\gamma_1,1))\vert_{1/D}(\gamma_2,1)=\eta^{2/D}\vert_{1/D}(\gamma_1\gamma_2,\sigma(\gamma_1,\sigma_2))$. Taking into account \eqref{eq:etaPowerTrans}, we obtain the desired identity.
\end{proof}

Now Proposition \ref{prop:PsiMultiplication} can be proved by letting $D$ tend to $+\infty$ in the above identity and using the formula for $\sigma(\gamma_1,\gamma_2)$ given above. (Notice that $\sigma(\gamma_1,\gamma_2)$ depends on $D$ which does not appear in the notation.)

\begin{coro}
\label{coro:charSlZ}
The map $\slZ\rightarrow\numZ/12\numZ$ that sends $\gamma$ to $\Psi(\gamma)+12\numZ$ is a surjective homomorphism. Moreover, the map $\chi_t\colon\slZ\rightarrow S^1$ defined by $\chi_t(\gamma)=\etp{\frac{t}{12}\Psi(\gamma)}$ is a linear character where $t\in\numZ$ and $\chi_t$ with $t=0,1,\dots,11$ are all different and constitute all linear characters of $\slZ$.
\end{coro}
\begin{proof}
The fact that the first map is a group homomorphism follows immediately from Proposition \ref{prop:PsiMultiplication}. It is surjective since $\Psi\tbtmat{1}{1}{0}{1}=1+12\numZ$. The second map is a linear character since the first map is a homomorphism. To see these 12 characters are all, note that $\slZ$ has a presentation with generators $T$, $S$ and relations $S^4=I$ and $S^2=(ST)^3$. Hence $\slZ$ has at most 12 linear characters. Since the characters $\gamma\mapsto\etp{\frac{t}{12}\Psi(\gamma)}$ are different (their values at $T$ are different), they are indeed all.
\end{proof}

\section{Explicit homomorphisms from $\Gamma_0(N)$ to $\numZ$}
\label{sec:HomGamma0N}
We define a map
\begin{align*}
\sigma_{N,l} \colon \Gamma_0(N) &\rightarrow \numZ \\
\gamma &\mapsto \Psi(\gamma)-\Psi(B_l^{-1}\gamma B_l),
\end{align*}
where $N$, $l$ are positive integers with $l\mid N$.

\begin{prop}
\label{prop:homoGamma0N}
The map $\sigma_{N,l}$ is a group homomorphism, and is nontrivial if $l>1$.
\end{prop}
\begin{proof}
Let $\gamma_i=\tbtmat{a_i}{b_i}{c_i}{d_i}$ be in $\Gamma_0(N)$ with $i=1,2$. Then $B_l^{-1}\gamma_i B_l=\tbtmat{a_i}{lb_i}{c_i/l}{d_i}$. Therefore by Proposition \ref{prop:PsiMultiplication} $\Psi(\gamma_1\gamma_2)=\Psi(\gamma_1)+\Psi(\gamma_2)+12$ if and only if $\Psi(B_l^{-1}\gamma_1 B_l\cdot B_l^{-1}\gamma_2 B_l)=\Psi(B_l^{-1}\gamma_1 B_l)+\Psi(B_l^{-1}\gamma_2 B_l)+12$. Similar equivalence holds if $12$ is replaced by $-12$, or by $0$. Hence $\sigma_{N,l}(\gamma_1\gamma_2)=\sigma_{N,l}(\gamma_1)+\sigma_{N,l}(\gamma_2)$, which shows $\sigma_{N,l}$ is a group homomorphism. To see it is nontrivial when $l>1$, note that $\sigma_{N,l}(T)=1-l$.
\end{proof}

This is the primary fact about $\sigma_{N,l}$ we need in the proof of the main theorem. However, $\sigma_{N,l}$ itself is an interesting object to study. For instance, it provides a method to find and prove identities concerning Dedekind sums. We will get back to this problem in the final section.

\section{Some preliminaries on group presentations}
\label{sec:groupPresentation}
We recall some elements in the theory of presentations of Fuchsian groups. The following lemma is an elementary fact about central extensions, which allows us to translate presentations of transformation groups into that of matrix groups.
\begin{lemm}
Let $G$ be a group and $g$ be an element of order $n\in\numgeq{Z}{1}$ in the center of $G$. Let $F$ be the free group with a basis $X$ and $f\colon F\rightarrow G/\langle g\rangle$ be a surjective group homomorphism. Let $R$ be a subset of $\ker f$ such that the normal closure of $R$ is $\ker f$. Let $x_g$ be a letter not in $X \cup X^{-1}$ and $F'$ be the free group with basis $X \cup \{x_g\}$. Choose a homomorphism $\tau\colon F'\rightarrow G$ as follows: for $x\in X$, $\tau(x)$ is any element in the coset $f(x)$ and $\tau(x_g)=g$. Set
\begin{align*}
R'&=\{r\cdot x_g^{-i}\colon r\in R,\, \tau(r)=g^i\,(0\leq i < n)\},\\
R''&=\{x_gxx_g^{-1}x^{-1}\colon x\in X\}.
\end{align*}
Then $\tau$ is surjective and the normal closure of $R'\cup R''\cup \{x_g^n\}$ is $\ker \tau$.
\end{lemm}
In another words, if $G/\langle g\rangle$ has a presentation $(X;R)$, then $G$ has a presentation $(X\cup\{x_g\};R'\cup R''\cup \{x_g^n\})$.
\begin{proof}
Let $K$ be the normal closure of $R'\cup R''\cup \{x_g^n\}$ in $F'$. The facts that $\tau$ is surjective and $K\subseteq \ker\tau$ are immediate. It remains to prove $\ker\tau \subseteq K$. Let $w \in \ker\tau$ be arbitrary and write it as a reduced word $w=y_1y_2\dots y_k$ with $y_i\in X\cup X^{-1}\cup\{x_g,x_g^{-1}\}$ and $y_iy_{i+1} \neq 1$. Since $x_gxx_g^{-1}x^{-1} \in K$ we have
\begin{equation*}
\prod_{i=1}^{k}y_iK=\left(\prod_{i=1}^{k'}y_i'K\right)x_g^tK
\end{equation*}
in the quotient group $F'/K$, where $y_i'\in X\cup X^{-1}$ and $t\in\numZ$. Since $K\subseteq \ker\tau$ we have
\begin{equation*}
\tau(y_1'y_2'\dots y_{k'}'x_g^t)=\tau(y_1y_2\dots y_k)=1.
\end{equation*}
Therefore $f(y_1'y_2'\dots y_{k'}')=1\cdot\langle g\rangle$, which means $y_1'y_2'\dots y_{k'}'$ is in the normal closure of $R$ in $F$. It follows that $y_1'y_2'\dots y_{k'}'x_g^tK= x_g^{t'}K$ for some $t'\in\numZ$ since $K$ contains $R'$ and $R''$. Now we have $g^{t'}=\tau(x_g^{t'})=\tau(y_1'y_2'\dots y_{k'}'x_g^t)=1$ so $n \mid t'$. Thus, $K=x_g^{t'}K= y_1'y_2'\dots y_{k'}'x_g^tK = wK$, that is, $w\in K$.
\end{proof}
If $G$ is a subgroup of $\slZ$, then its image under the natural projection $\slZ\rightarrow\pslZ$ will be denoted by $\overline{G}$, and $\overline{g}=g\cdot\{\pm I\}$ for $g \in \slZ$.
\begin{coro}
\label{coro:presMatrixGrp}
Let $G$ be a finite-index subgroup of $\slZ$ and suppose that $\{\overline{g_1},\dots,\overline{g_r},\overline{h_1},\dots \overline{h_k}\}$ is an independent set of generators of $\overline{G}$ (that is, $\overline{G}$ is a free product of each $\langle \overline{g_i}\rangle$, $\langle \overline{h_i}\rangle$) with relations $\overline{h_i}^{m_i}=\overline{I}$ ($m_i=2$ or $3$). Such set of generators always exists by Kurosh subgroup theorem. Then after a possible replacement of $h_i$ with $-h_i$, $G$ has a presentation with generators $g_1,\dots,g_r,h_1,\dots,h_k,-I$ and relations $h_i^{m_i}=-I$, $(-I)g_i=g_i(-I)$, $(-I)h_i=h_i(-I)$, $(-I)^2=I$.
\end{coro}
\begin{proof}
An application of the previous lemma with $G=\slZ$ and $g=-I$.
\end{proof}
\begin{rema}
\label{rema:presMatrixGrp}
In the following, we will always use the independent set of generators for $\overline{\Gamma_0(N)}$ given by Kulkarni's algorithm \cite[Theorem 13.2 and Theorem 6.1]{Kul91}. The corresponding set of generators of $\Gamma_0(N)$ can be obtained by the SageMath \cite{Sage} code \lstinline{Gamma0(N).gens()}.
\end{rema}

Let $G=\ast_{i \in I}G_i$ be a free product, $[G,G]$ be the commutator subgroup and $\mathop{\mathrm{Ab}}(G)=G/[G,G]$ be the abelianization. Then $\mathop{\mathrm{Ab}}(G)$ is isomorphic to the direct sum of the groups $\mathop{\mathrm{Ab}}(G_i),\,i\in I$ written additively. This can be verified directly by definitions. Moreover, a finite-index subgroup of a finitely presented group is finitely presented (c.f. \cite[Proposition 4.2, Chapter II]{LS77}). Thus, if $\overline{G}$ is a finite-index subgroup of $\pslZ$ which (by Kurosh subgroup theorem) is written as a free product of $r$ copies of $\numZ$, $e_2$ copies of $\numZ/2\numZ$ and $e_3$ copies of $\numZ/3\numZ$, then $r$, $e_2$, $e_3$ are finite and by abelianizing they are independent of the specific decomposition. The \emph{measure} of $\overline{G}$ is defined by
\begin{equation*}
\mu(\overline{G})=(r+e_2+e_3)-1-\frac{e_2}{2}-\frac{e_3}{3}.
\end{equation*}
This concept has been used since the time of Poincar\'{e} and Klein. The following theorem concerning measures is a special case of \cite[Theorem 3]{HKS71}.
\begin{thm}
\label{thm:RH}
Let $G$ be a finite-index subgroup of $\slZ$. Then $\mu(\overline{G})=[\pslZ\colon \overline{G}]\cdot \mu(\pslZ)$.
\end{thm}
The expression $[\pslZ\colon \overline{G}]$ means the index of $\overline{G}$ in $\pslZ$. For $\Gamma_0(N)$ we have (c.f. \cite[Corollary 6.2.13]{CS17})
\begin{equation*}
[\pslZ\colon \overline{\Gamma_0(N)}]=[\slZ\colon \Gamma_0(N)]=N\prod_{p \mid N}\left(1+\frac{1}{p}\right),
\end{equation*}
where $p$ denotes primes. Theorem \ref{thm:RH} can be derived from the Riemann-Hurwitz formula applying to the (holomorphic) branched covering $\overline{G}\backslash\uhp^{\ast}\rightarrow\pslZ\backslash\uhp^{\ast}$, while Hoare, Karrass and Solitar used purely combinatorial methods. Writing out explicitly, we have
\begin{equation}
\label{eq:RH}
r=\frac{1}{6}N\prod_{p \mid N}\left(1+\frac{1}{p}\right)+1-\frac{1}{2}e_2-\frac{2}{3}e_3,
\end{equation}
where $r$, $e_2$, $e_3$ are the numbers of factors $\numZ$, $\numZ/2\numZ$, $\numZ/3\numZ$ in the free-product decomposition of $\overline{\Gamma_0(N)}$ respectively.

\section{Explicit linear unitary characters of $\Gamma_0(N)$}
\label{sec:CharGamma0N}
For a group $G$, let $\widehat{G}$ be the group of its linear unitary characters over $\numC$, with the group composition being pointwise multiplication. Thus, $\widehat{\numZ/N\numZ^\times}$ is isomorphic to the group of Dirichlet characters modulo $N$ and $\widehat{\Gamma_0(N)}$ is what we are studying.

We separate the main theorem into two parts. The former is the relatively simpler part, which is slightly stronger than the one stated in Introduction.
\begin{thm}
\label{thm:main1}
Let $N$ be a positive integer and $t$ be the number of positive divisors of $N$. Then the map
\begin{align}
\label{eq:allGamma0Nchar}
\widehat{\numZ/N\numZ^\times}\times \numZ/12\numZ \times \numR^{t-1} &\rightarrow \widehat{\Gamma_0(N)}\\
(\chi, r_1, (r_l)_{1<l \mid N}) &\mapsto \left(\gamma\mapsto\chi(d)\etp{\frac{r_1}{12}\Psi(\gamma)}\etp{\sum_{1<l\mid N}r_l\sigma_{N,l}(\gamma)}\right)\notag
\end{align}
is a group homomorphism, where $\gamma=\tbtmat{*}{*}{*}{d}$.
\end{thm}
\begin{proof}
Let $(\chi, r_1, (r_l)_{1<l \mid N}) \in \widehat{\numZ/N\numZ^\times}\times \numZ/12\numZ \times \numR^{t-1}$ and $\gamma_1,\gamma_2 \in \Gamma_0(N)$ be arbitrary. Write $\gamma_i=\tbtmat{a_i}{b_i}{c_i}{d_i}$ and $\chi(\gamma_i)=\chi(d_i)$, $i=1,2$. Then
\begin{equation*}
\chi(\gamma_1\gamma_2)=\chi(c_1b_2+d_1d_2)=\chi(d_1)\chi(d_2)=\chi(\gamma_1)\chi(\gamma_2).
\end{equation*}
Therefore the map $\gamma\mapsto\chi(\gamma)=\chi(d)$ is in $\widehat{\Gamma_0(N)}$. By Corollary \ref{coro:charSlZ} and Proposition \ref{prop:homoGamma0N}, the maps $\gamma\mapsto\etp{\frac{r_1}{12}\Psi(\gamma)}$ and $\gamma\mapsto\etp{\sum_{1<l\mid N}r_l\sigma_{N,l}(\gamma)}$ are also in $\widehat{\Gamma_0(N)}$. The image of $(\chi, r_1, (r_l)_{1<l \mid N})$ is the product of these three characters, and hence is in $\widehat{\Gamma_0(N)}$. The fact that \eqref{eq:allGamma0Nchar} is a homomorphism is immediate.
\end{proof}

The latter part of the main theorem can be restated as
\begin{thm}
\label{thm:main2}
The homomorphism \eqref{eq:allGamma0Nchar} is surjective if and only if $N=1,2,3,4,5,6,7,8,10,12,13$.
\end{thm}
\begin{proof}
The ``if'' part. If $N=1$, the assertion follows from Corollary \ref{coro:charSlZ}. Assume that $N\in\{2,3,4,5,6,7,8,10,12,13\}$. Using Kulkarni's algorithm \cite[Theorem 13.2 and Theorem 6.1]{Kul91} we obtain an independent set of generators of $\overline{\Gamma_0(N)}$ and then using Corollary \ref{coro:presMatrixGrp} we obtain a presentation of $\Gamma_0(N)$. This can be done by the SageMath \cite{Sage} code \lstinline{Gamma0(N).gens()}. The generators obtained in this way are listed in Table \ref{table:gens} which are grouped by their orders.

\begin{table}[ht]
\centering
\caption{For each $N$ a set of generators of $\Gamma_0(N)$, where the column titled $\infty$ contains generators of infinite order, that titled $e$ ($2$ or $3$) contains generators $\gamma$ such that $\gamma^e=-I$, and the generator $-I$ is excluded\label{table:gens}}
\begin{tabular}{llll}
\toprule
$N$ & $\infty$ & $2$ & $3$ \\
\midrule
$2$ & $\tbtmat{1}{1}{0}{1}$ & $\tbtmat{1}{-1}{2}{-1}$ & $\emptyset$ \\
$3$ & $\tbtmat{1}{1}{0}{1}$ & $\emptyset$ & $\tbtmat{-1}{1}{-3}{2}$ \\
$4$ & $\tbtmat{1}{1}{0}{1}$,$\tbtmat{3}{-1}{4}{-1}$& $\emptyset$ & $\emptyset$ \\
$5$ & $\tbtmat{1}{1}{0}{1}$ & $\tbtmat{2}{-1}{5}{-2}$, $\tbtmat{3}{-2}{5}{-3}$ & $\emptyset$ \\
$6$ & $\tbtmat{1}{1}{0}{1}$,$\tbtmat{5}{-1}{6}{-1}$,$\tbtmat{7}{-3}{12}{-5}$& $\emptyset$ & $\emptyset$ \\
$7$ & $\tbtmat{1}{1}{0}{1}$ & $\emptyset$ & $\tbtmat{-2}{1}{-7}{3}$, $\tbtmat{-4}{3}{-7}{5}$ \\
$8$ & $\tbtmat{1}{1}{0}{1}$,$\tbtmat{5}{-1}{16}{-3}$,$\tbtmat{5}{-2}{8}{-3}$& $\emptyset$ & $\emptyset$ \\
$10$ & $\tbtmat{1}{1}{0}{1}$,$\tbtmat{19}{-7}{30}{-11}$,$\tbtmat{11}{-5}{20}{-9}$ & $\tbtmat{3}{-1}{10}{-3}$,$\tbtmat{7}{-5}{10}{-7}$ & \\
$12$ & $\tbtmat{1}{1}{0}{1}$,$\tbtmat{7}{-1}{36}{-5}$,$\tbtmat{19}{-4}{24}{-5}$,$\tbtmat{17}{-5}{24}{-7}$,$\tbtmat{7}{-3}{12}{-5}$ & $\emptyset$ & $\emptyset$ \\
$13$ & $\tbtmat{1}{1}{0}{1}$ & $\tbtmat{5}{-2}{13}{-5}$,$\tbtmat{8}{-5}{13}{-8}$ & $\tbtmat{-3}{1}{-13}{4}$,$\tbtmat{-9}{7}{-13}{10}$ \\
\bottomrule
\end{tabular}
\end{table}

Now any character $v\in\widehat{\Gamma_0(N)}$ is uniquely determined by assigning a value in $S^1$ to each generator of infinite order, a value in $\{\pm 1\}$ to $-I$, a value in $\{\pm 1,\pm \rmi\}$ to each generator $\gamma$ such that $\gamma^e=-I$ with $e=2$ and a value in $\{\etp{n/6}\colon n=0,1,2,3,4,5\}$ to each generator $\gamma$ such that $\gamma^e=-I$ with $e=3$, subject to the conditions $v(\gamma)^e=v(-I)$. For each $N$ in the table and each $v\in\widehat{\Gamma_0(N)}$ we now construct $(\chi, r_1, (r_l)_{1<l \mid N}) \in \widehat{\numZ/N\numZ^\times}\times \numZ/12\numZ \times \numR^{t-1}$ whose image under \eqref{eq:allGamma0Nchar} is $v$. Note that for any $N$ in the table, if we list the generators of order $\infty$ as $g_1,\ldots,g_r$, then $r=t-1$ and the square matrix $\left(\sigma_{N,l}(g_i)\right)_{1\leq i \leq r,1<l\mid N}$ is nonsingular, which can be verified by a direct calculation (using any computer algebra system).

For $N=4,6,8,12$, let $\chi$ be the trivial character in $\widehat{\numZ/N\numZ^\times}$, and let $r_1=0$ or $1$ according to $v(-I)=1$ or $-1$. Since $\left(\sigma_{N,l}(g_i)\right)_{1\leq i \leq r,1<l\mid N}$ is nonsingular, there exists $(r_l)_{1<l \mid N}$ such that the image of $(\chi, r_1, (r_l)_{1<l \mid N})$ is $v$.

For $N=2,3$, set $\gamma_0=\tbtmat{1}{-1}{2}{-1}$ if $N=2$, or $\tbtmat{-1}{1}{-3}{2}$ if $N=3$. Let $\chi$ be the trivial character in $\widehat{\numZ/N\numZ^\times}$. If $N=2$ and $v(\gamma_0)=\etp{m/4}$, $m=0,1,2,3$, then let $r_1=4-m$; if $N=3$ and $v(\gamma_0)=\etp{m/6}$, $m=0,1,2,3,4,5$, then let $r_1=m$. Choose $r_N \in \numR$ such that $\etp{\frac{r_1}{12}\Psi\tbtmat{1}{1}{0}{1}}\etp{r_N\sigma_{N,N}\tbtmat{1}{1}{0}{1}}=v\tbtmat{1}{1}{0}{1}$, then the image of $(\chi, r_1, r_N)$ is $v$.

The constructions for $N=5,7$ are similar, so we omit the case $N=5$ and consider $N=7$. Set $z_1=v\tbtmat{-2}{1}{-7}{3}$ and $z_2=v\tbtmat{-4}{3}{-7}{5}$. Then $z_1^3=z_2^3=1$ (when $v(-I)=1$) or $z_1^3=z_2^3=-1$ (when $v(-I)=-1$). Set $z_1/z_2=\etp{j/3}$ with $j=0,1$ or $2$. Note that $\Psi\tbtmat{-2}{1}{-7}{3}=\Psi\tbtmat{-4}{3}{-7}{5}=2$. Let $\chi$ be the character in $\widehat{\numZ/7\numZ^\times}$ that maps $3+7\numZ$ to $\etp{j/3}$. Then it maps $5+7\numZ$ to $\etp{2j/3}$. Therefore, if we choose $r_1$ such that $\chi\tbtmat{-2}{1}{-7}{3}\etp{\frac{r_1}{12}\Psi\tbtmat{-2}{1}{-7}{3}}=z_1$, then automatically $\chi\tbtmat{-4}{3}{-7}{5}\etp{\frac{r_1}{12}\Psi\tbtmat{-4}{3}{-7}{5}}=z_2$. Now let $r_7 \in \numR$ satisfy $\chi\tbtmat{1}{1}{0}{1}\etp{\frac{r_1}{12}\Psi\tbtmat{1}{1}{0}{1}+r_7\sigma_{7,7}\tbtmat{1}{1}{0}{1}}=v\tbtmat{1}{1}{0}{1}$. It follows that the image of $(\chi, r_1, r_7)$ is $v$ since $\sigma_{7,7}\tbtmat{-2}{1}{-7}{3}=\sigma_{7,7}\tbtmat{-4}{3}{-7}{5}=0$.

The proof for the case $N=10$ is also similar to the case $N=5$ so we omit it. One uses the fact that $\left(\sigma_{10,l}(g_i)\right)_{1\leq i \leq 3,1<l\mid 10}$ is nonsingular, uses the character $\chi\in\widehat{\numZ/10\numZ^\times}$ that maps $1,3,7,9$ to $1,\rmi,-\rmi,-1$ respectively, and chooses $r_1\in\{0,1,2,3\}$.

For $N=13$, set\footnote{In this case, the notations $h_i$,$g_i$ have different meanings than Corollary \ref{coro:presMatrixGrp}.} $T=\tbtmat{1}{1}{0}{1}$, $g_1=\tbtmat{5}{-2}{13}{-5}$, $g_2=\tbtmat{8}{-5}{13}{-8}$, $h_1=\tbtmat{-3}{1}{-13}{4}$ and $h_2=\tbtmat{-9}{7}{-13}{10}$. Set $w_i=v(g_i)$ and $z_i=v(h_i)$ with $i=1,2$. Then $w_1^2=w_2^2=z_1^3=z_2^3=1$ (when $v(-I)=1$) or $w_1^2=w_2^2=z_1^3=z_2^3=-1$ (when $v(-I)=-1$). For any $x,y\in\numZ/12\numZ$ there exists a unique $\chi\in\widehat{\numZ/13\numZ^\times}$ such that $\chi(2+13\numZ)=\etp{x/12}$, and consequently letting $r_1=y$ we have
\begin{align*}
\chi(g_1)\etp{\frac{r_1}{12}\Psi(g_1)}&=\etp{\frac{x-y}{4}}\\
\chi(g_2)\etp{\frac{r_1}{12}\Psi(g_2)}&=\etp{\frac{-x-y}{4}}\\
\chi(h_1)\etp{\frac{r_1}{12}\Psi(h_1)}&=\etp{\frac{x+y}{6}}\\
\chi(h_2)\etp{\frac{r_1}{12}\Psi(h_2)}&=\etp{\frac{-x+y}{6}}.
\end{align*}
Therefore, for any $(w_1,w_2,z_1,z_2)$ we can find $(x,y)$ and hence $(\chi,r_1)$ such that $\chi(g_i)\etp{\frac{r_1}{12}\Psi(g_i)}=w_i$ and $\chi(h_i)\etp{\frac{r_1}{12}\Psi(h_i)}=z_i$ with $i=1,2$. Let $r_{13} \in \numR$ satisfy $\chi(T)\etp{\frac{r_1}{12}\Psi(T)+r_{13}\sigma_{13,13}(T)}=v(T)$. Then the image of $(\chi, r_1, r_{13})$ is $v$ since $\sigma_{13,13}(h_i)=\sigma_{13,13}(g_i)=0$ for $i=1,2$.

This concludes the proof of the ``if'' part, so we turn to the ``only if'' part. Let $N$ be any positive integer; we shall prove that if $N \neq 1,2,3,4,5,6,7,8,10,12,13$, then \eqref{eq:allGamma0Nchar} is not surjective. We use the notations of \eqref{eq:RH}, Corollary \ref{coro:presMatrixGrp} (with $G=\Gamma_0(N)$) and Remark \ref{rema:presMatrixGrp} and break the argument into two cases:
\begin{equation*}
\text{(a)}\, 2^{e_2+1}3^{e_3}>12N\prod_{p \mid N}\left(1-\frac{1}{p}\right),\qquad \text{(b)}\, 2^{e_2+1}3^{e_3}\leq 12N\prod_{p \mid N}\left(1-\frac{1}{p}\right).
\end{equation*}

For case (a), note that the set $\widehat{\numZ/N\numZ^\times}\times \numZ/12\numZ$ has cardinality $12N\prod_{p \mid N}\left(1-\frac{1}{p}\right)$, and that $\sigma_{N,l}(h_i)=0$ (for the notation $h_i$, see Corollary \ref{coro:presMatrixGrp}). Therefore the images of all $(\chi, r_1, (r_l)_{1<l \mid N})$ give at most $12N\prod_{p \mid N}\left(1-\frac{1}{p}\right)$ characters on the subgroup generated by $\{h_i\}$. However, this group has exactly $2^{e_2+1}3^{e_3}$ characters according to the presentation in Corollary \ref{coro:presMatrixGrp} and Remark \ref{rema:presMatrixGrp}. It follows that \eqref{eq:allGamma0Nchar} is not surjective.

For case (b), we have $e_2 \leq \log_2(12N\prod(1-1/p))-1$ and $e_3 \leq \log_3(12N\prod(1-1/p))$. Using these inequalities and \eqref{eq:RH} we obtain
\begin{multline*}
r \geq \frac{1}{6}N\prod_{p \mid N}\left(1+\frac{1}{p}\right)-\left(\frac{1}{2\log2}+\frac{2}{3\log3}\right)\log{\left(N\prod_{p \mid N}\left(1-\frac{1}{p}\right)\right)}
-\left(\frac{1}{6}+\frac{\log 3}{2\log 2}+\frac{4\log 2}{3\log 3}\right).
\end{multline*}
Rounding the constants and omitting the products over $p$ we find that $r>\frac{1}{6}N-1.4\log N-2$. It is elementary that $t$, the number of positive divisors of $N$, is less than $2\sqrt{N}$. Consequently, by elementary calculus, if $N \geq 237$, then
\begin{equation*}
r > \frac{1}{6}N-1.4\log N-2 > 2\sqrt{N} - 1 > t-1.
\end{equation*}
If $14 \leq N \leq 236$ or $N=9,11$, we also have $r > t-1$ by a direct calculation (using SageMath). For any fixed $(\chi, r_1)\in\widehat{\numZ/N\numZ^\times}\times \numZ/12\numZ$, consider the map
\begin{align*}
f_{\chi,r_1}\colon \numR^{t-1} &\rightarrow (S^1)^r\\
(r_l)_{1<l \mid N} &\mapsto \left(z_1,z_2,\dots,z_r\right),
\end{align*}
where (the notation $g_j$ is introduced in Corollary \ref{coro:presMatrixGrp})
\begin{equation*}
z_j=\chi(g_j)\etp{\frac{r_1}{12}\Psi(g_j)}\etp{\sum\nolimits_{1<l\mid N}r_l\sigma_{N,l}(g_j)}.
\end{equation*}
It is a smooth map with a constant rank equal to the rank of $\left(\sigma_{N,l}(g_j)\right)_{1\leq j \leq r,1<l\mid N}$. Since $r>t-1$, it can not be a submersion. Thus\footnote{One can modify the argument in the proof of \cite[Theorem 4.14(a)]{Lee13} slightly to obtain a detailed proof of this.} $f_{\chi,r_1}(\numR^{t-1})$ is a countable union of nowhere dense sets in $(S^1)^r$. Therefore we have
\begin{equation*}
\bigcup_{(\chi,r_1)\in\widehat{\numZ/N\numZ^\times}\times \numZ/12\numZ}f_{\chi,r_1}(\numR^{t-1}) \neq (S^1)^r.
\end{equation*}
To see this, suppose the contrary; taking the closure and then the complement of both sides we find that the intersection of a countable collection of dense open sets is empty, which contradicts the Baire category theorem applied to the locally compact Hausdorff space $(S^1)^r$. It follows immediately that \eqref{eq:allGamma0Nchar} is not surjective, which concludes the proof of the ``only if'' part.
\end{proof}

\section{Related results and Open questions}
\label{sec:OpenQues}
The major tool we introduced to prove the main theorem (i.e. Theorems \ref{thm:main1} and \ref{thm:main2}) is the homomorphism $\sigma_{N,l}$ defined in Section \ref{sec:HomGamma0N}. This homomorphism has many other applications. In this section, we examine some of these together with some open questions and conjectures which are not related to proving the main theorem.
\subsection{The kernel of $\sigma_{N,l}$}
The homomorphism $\sigma_{N,l}$ can be used to derive identities on Dedekind sums. This is based on the following observation: for any $\gamma=\tbtmat{a}{b}{c}{d}\in\Gamma_0(N)$ with $c>0$, $\gamma\in\ker\sigma_{N,l}$ if and only if
\begin{equation*}
12s(-d,c)=12s(-d,c/l)+\frac{a+d}{c}(l-1).
\end{equation*}
As a simple illustration, take a $\gamma\in\Gamma_0(N)$ such that $\gamma^3=\pm I$; then $\gamma\in\ker\sigma_{N,l}$. In this way we can find that $s(d,c)=\frac{c-1}{12c}$ if $d^2+d+1 \equiv 0 \bmod c$ and $c>0$, which is a special case of \cite[Sat 13]{Rad56}. Rademacher used his original $\Psi$ function and the theory of quadratic forms, while our method is different from Rademacher's.

In view of the above consideration, it is useful to describe $\ker\sigma_{N,l}$ explicitly. Note that $\ker\sigma_{N,l}=\ker\sigma_{l,l}\cap\Gamma_0(N)$, so we only consider $\ker\sigma_{N,N}$ without loss of generality. Let us use the notations and presentations of $\Gamma_0(N)$ in Corollary \ref{coro:presMatrixGrp} and Remark \ref{rema:presMatrixGrp} and first explain some general principles. For any $\gamma\in\Gamma_0(N)$ with $\gamma \neq \pm I$, by the normal form theorem\footnote{Some authors use this as the definition of free products, c.f. \cite{Mac58}} for free products over arbitrary index set (an immediate generalization of \cite[Theorem 1.2, Chapter IV]{LS77} which deals with free products of two factors), we have $\overline{\gamma}=\prod_{i=1}^n\overline{\gamma_i}^{\alpha_i}$ where $n \geq 1$, $\gamma_i\in\{g_1,\dots,g_r,h_1,\dots,h_k\}$, $\gamma_i\neq \gamma_{i+1}$ and $\alpha_i\in\numZ-\{0\}$ if $\gamma_i\in\{g_1,\dots,g_r\}$; $\alpha_i=1$ if $\gamma_i$ equals some $h_j$ with $m_j=2$; $\alpha_i\in\{1,2\}$ if $\gamma_i$ equals some $h_j$ with $m_j=3$. Moreover, such expression is unique. It follows that $\gamma=(\varepsilon I)\cdot\prod_{i=1}^n\gamma_i^{\alpha_i}$ with $\varepsilon=\pm 1$ and such decomposition is also unique. We call this the normal form of $\gamma\in\Gamma_0(N)$ with respect to the given presentation in Corollary \ref{coro:presMatrixGrp} and Remark \ref{rema:presMatrixGrp}. Since $\sigma_{N,N}(h_i)=0$ and $\sigma_{N,N}(\pm I)=0$, we conclude that $\gamma\in\ker\sigma_{N,N}$ if and only if
\begin{equation*}
\sum_{\twoscript{1\leq i \leq n}{\gamma_i\in\{g_1,\dots,g_r\}}}\alpha_i\sigma_{N,N}(\gamma_i)=0.
\end{equation*}

In some special cases, this leads to clear description of $\ker\sigma_{N,N}$.
\begin{prop}
\label{prop:kersigmaNN}
Let $N$ be a positive integer and use the presentation given in Corollary \ref{coro:presMatrixGrp} (with $G=\Gamma_0(N)$) and Remark \ref{rema:presMatrixGrp}. If there is exactly one $1\leq j \leq r$ such that $\sigma_{N,N}(g_j)\neq 0$, then
\begin{equation*}
\ker\sigma_{N,N}=\left\{\text{normal form }(\varepsilon I)\cdot\prod\nolimits_{i=1}^n\gamma_i^{\alpha_i}\in\Gamma_0(N)\colon \sum\nolimits_{1\leq i\leq n,\,\gamma_i=g_j}\alpha_i=0\right\}.
\end{equation*}
In particular, the normal closure of the set of generators excluding $g_j$ is contained in $\ker\sigma_{N,N}$.
\end{prop}

Among $2 \leq N < 1000$, we find out all $N$ such that the condition in the preceding proposition holds by a SageMath program:
\begin{equation*}
N=2,3,4,5,7,9,13,25.
\end{equation*}

We now give a generalization of a property \cite[Theorem 3.11]{Apo76} of Dedekind sums which plays an important role in the theory of modular functions.
\begin{prop}
\label{prop:dedekindsumIdentity}
Let $N \in \{2,3,4,5,7,9,13,25\}$, $c,d\in\numZ$ with $N \mid c$, $\gcd(c,d)=1$ and $c>0$. Let $a$ be an integer such that $ad \equiv 1\bmod c$. Then
\begin{equation}
\label{eq:dedekindsumIdentity}
\left(\frac{a+d}{c}-12s(d,c)\right)-\left(\frac{a+d}{c/N}-12s(d,c/N)\right)\in(N-1)\numZ.
\end{equation}
\end{prop}
\begin{proof}
Set $b=(ad-1)/c$ and $\gamma=\tbtmat{a}{b}{c}{d}$. Then $\gamma\in\Gamma_0(N)$. Since $N \in \{2,3,4,5,7,9,13,25\}$ it satisfies the condition in Proposition \ref{prop:kersigmaNN} and $g_j=T=\tbtmat{1}{1}{0}{1}$. Therefore if in the normal form of $\gamma$ the sum of the exponents of all $T$ is $-n\in\numZ$, then $T^n\gamma\in\ker\sigma_{N,N}$, which means
\begin{equation*}
\frac{a+nc+d}{c}+12s(-d,c)-3=\frac{a+nc+d}{c/N}+12s(-d,c/N)-3.
\end{equation*}
A simplification shows that
\begin{equation*}
\left(\frac{a+d}{c}-12s(d,c)\right)-\left(\frac{a+d}{c/N}-12s(d,c/N)\right)=(N-1)n\in(N-1)\numZ.
\end{equation*}
\end{proof}
\begin{rema}
In \cite[Theorem 3.11]{Apo76}, Apostol obtained this relation only for $N=3,5,7,13$. Moreover, our proof is shorter and more elegant.
\end{rema}
\begin{rema}
Proposition \ref{prop:dedekindsumIdentity} can be used to prove certain eta-quotients are modular functions. Precisely, $$\left(\frac{\eta(N\tau)}{\eta(\tau)}\right)^{\frac{24}{N-1}}$$ is a modular function under $\Gamma_0(N)$ whenever $N-1$ is a positive divisor of $24$. The proof for $N=2,3,5,7,13$ can be found in \cite[Theorem 4.9]{Apo76} and that for other $N$ is similar. This is a special case of Newman's theorem \cite{New59}.
\end{rema}

\subsection{The image of $\sigma_{N,l}$}
Keep the notations of Corollary \ref{coro:presMatrixGrp} and Remark \ref{rema:presMatrixGrp} with $G=\Gamma_0(N)$ and assume that $N>1$ and $1<l \mid N$. Then there is a unique positive integer $\beta=\beta(N,l)$ such that $\mathrm{im}(\sigma_{N,l})=\beta\numZ$, where $\mathrm{im}(f)$ means the image of $f$. It is immediate that $\beta=\gcd(\sigma_{N,l}(g_1),\ldots,\sigma_{N,l}(g_r))$ since $\sigma_{N,l}(h_j)=0$ for any $1 \leq j \leq k$ and $\sigma_{N,l}(-I)=0$. Thus, it is easy to algorithmically compute $\beta(N,l)$ for any $N,l$. Moreover, there always exists some $\gamma_0\in\Gamma_0(N)$ such that $\sigma_{N,l}(\gamma_0)=\beta$. This $\gamma_0$ will be needed if one considers, for instance, explicitly constructing noncongruence subgroups or explicitly computing Fourier coefficients of Eisenstein series on these noncongruence subgroups.

After some experiments we raise some conjectures on $\beta(N,l)$:
\begin{conj}
\label{conj1}
Let $N \geq 2$ and $1 < l \mid N$. Then $\beta(N,l)=\beta(l,l)$.
\end{conj}
This has been verified to be true for $l,N<1000$ by a SageMath program.

If we acknowledge the above conjecture, then we need only to study $\beta(N,N)$. By abuse of notation, we write $\beta(N)=\beta(N,N)$. The following conjecture concerning $\beta(N)$ is striking to us:
\begin{conj}
\label{conj2}
Let $N \geq 2$ be an integer.
\begin{enumerate}
  \item If $N \not\equiv 1,9 \bmod{24}$, then $\beta(N)$ is a periodic function of $N$ with a period $24$. The values $\beta(N)$ are shown in Table \ref{table:betaN}.
  \item If $N \equiv 9 \bmod{24}$, then $\beta(N)=8$ if $N$ is a perfect square and $\beta(N)=4$ otherwise.
  \item If $N \equiv 1 \bmod{24}$, then $\beta(N)=24$ if $N$ is a perfect square and $\beta(N)=12$ otherwise. 
\end{enumerate}
\end{conj}
\begin{table}[ht]
\centering
\caption{The values $\beta(N)$ for $N\geq2$, where $N \equiv R \bmod 24$, $1\leq R \leq 24$ \label{table:betaN}}
\begin{tabular}{llllllll}
\toprule
$R$ & $\beta(N)$ & $R$ & $\beta(N)$ & $R$ & $\beta(N)$ & $R$ & $\beta(N)$ \\
\midrule
$1$ & $12$ or $24$  &  $7$ &  $6$         & $13$ &  $12$  &  $19$ & $6$  \\
$2$ & $1$           &  $8$ &  $1$         & $14$ &  $1$   &  $20$ & $1$ \\
$3$ & $2$           &  $9$ &  $4$ or $8$  & $15$ &  $2$   &  $21$ & $4$ \\
$4$ & $3$           &  $10$&  $3$         & $16$ &  $3$   &  $22$ & $3$ \\
$5$ & $4$           &  $11$&  $2$         & $17$ &  $4$   &  $23$ & $2$ \\
$6$ & $1$           &  $12$&  $1$         & $18$ &  $1$   &  $24$ & $1$ \\
\bottomrule
\end{tabular}
\end{table}

Conjecture \ref{conj2}(1) has been verified to be true for $N \leq 960$ and Conjectures \ref{conj2}(2) and (3) true for $N<3000$ by a SageMath program.

\subsection{The rank of $\left(\sigma_{N,l}(g_j)\right)_{1\leq j \leq r,1<l\mid N}$}
In the proof of Theorem \ref{thm:main2} we have used the $r\times(t-1)$ matrix $\left(\sigma_{N,l}(g_j)\right)_{1\leq j \leq r,1<l\mid N}$ several times. For $g_j$ and $r$ see Corollary \ref{coro:presMatrixGrp} and Remark \ref{rema:presMatrixGrp}; $t$ is the number of positive divisors of $N$. The rank of this matrix is related to the question how large would the image of \eqref{eq:allGamma0Nchar} be. Theorem \ref{thm:main2} gives a partial answer by finding out when the image equals the whole group of unitary characters. To answer this question for larger $N$, it is necessary to know the rank first, on which we have the following conjecture.
\begin{conj}
\label{conj3}
Let the notations be as above and suppose $N \geq 2$, then
\begin{equation*}
\mathop{\mathrm{rank}}{\left(\sigma_{N,l}(g_j)\right)_{1\leq j \leq r,1<l\mid N}}=t-1.
\end{equation*}
\end{conj}
This has been verified to be true for $N<1000$ by a SageMath program. Note that if $N$ is a prime, then it holds trivially.

\section*{Statements and Declarations}
The author declares that he has no conflicts of interest in the research presented in this manuscript.

All data and SageMath code supporting the findings of this study are available in \cite{Zhu25}.

\bibliographystyle{amsplain}
\bibliography{ref}
\end{document}